\documentclass{elsart3-1}

\usepackage{amsmath} 
\def\NN{\bf N}

\def\cP{\mathcal P}

\def\bw{{\bf w}}
\def\bx{{\bf x}}

\def\bv{{\bf v}}

\def\be{{\bf e}}
\def\br{{\bf r}}
\def\bg{{\bf g}}

\def\bgamma{\boldsymbol{\gamma}}

\def\varkappa{{\kappa}}

\def\bs{{\bf 0}}
\def\ad{{\rm and}}
\def\with{{\rm with}}
\def\where{{\rm where}}
\def\for{{\rm for}}

\usepackage[english,francais]{babel}

\newtheorem{e-proposition}[theorem]{Proposition}

\newtheorem{e-definition}[theorem]{Definition\rm}

\renewcommand{\theequation}{\arabic{equation}}
\def\og{\leavevmode\raise.3ex\hbox{$\scriptscriptstyle\langle\!\langle$~}}
\def\fg{\leavevmode\raise.3ex\hbox{~$\!\scriptscriptstyle\,\rangle\!\rangle$}}

\journal{the Acad\'emie des sciences}
\begin{document}
\centerline{}
\begin{frontmatter}




%
\selectlanguage{english}
\title{On the lower bound of the discrepancy  of $(t,s)$ sequences: I}



\author{Mordechay B. Levin},
\ead{mlevin@math.biu.ac.il}

\address{Department of Mathematics,
Bar-Ilan University, Ramat-Gan, 52900, Israel}



\medskip
\selectlanguage{francais}
\begin{center}
{\small Re\c{c}u le *****~; accept\'e apr\`es r\'evision le +++++\\
Pr\'esent\'e par £££££}
\end{center}

\selectlanguage{english}
\begin{abstract}
\selectlanguage{english}
\vskip 0.5\baselineskip
We find the exact lower bound of the discrepancy  of shifted   Niedereiter's  sequences.\\
{\it To cite this article: A. Nom1,  C. R.
Acad. Sci. Paris, Ser. I .}

\selectlanguage{francais}
\noindent{\bf R\'esum\'e} \vskip 0.5\baselineskip \noindent
{\bf Sur la limite inf\'erieure de la discr\'epance de $(t,s)$ s\'equences: I }\\
  Nous trouvons une limite inf\'erieure  de la  discr\'epance de    s\'equences d\'ecal\'es
 de Niedereiter. \\
{\it Pour citer cet article~: A. Nom1,  C. R.
Acad. Sci. Paris, Ser. I .}
\vskip 0.5\baselineskip
\end{abstract}
\end{frontmatter}

\selectlanguage{english}




\section{Introduction}
 Let $((\bx_{n})_{n \geq 1})$ be an  $s$-dimensional sequence in the unit cube $[0,1)^s$,
$J_{\bgamma}=[0,\gamma_1) \times \cdots \times [0,\gamma_s) \subseteq [0,1)^s $,
\begin{equation}\label{1}
\Delta( (\bx_{n})_{n=1}^{N}, J_{\bgamma}  )= \sum_{0 \leq n
  < N} {\bf 1} ( \bx_{n}, J_{\bgamma}) -N\gamma_1 \ldots \gamma_s,
\end{equation}
where ${\bf 1} ( \bx, J) =1$, if $ \bx \in J$ and ${\bf 1} (\bx, J) =0$, if $ \bx \notin J$. 
We define the star {\it discrepancy} of a 
$N$-point set $(\beta_{n})_{n=1}^{N}$ as
\begin{equation} \nonumber
  \emph{D}^{*}((\bx_{n})_{n=1}^{N}) = 
    \sup_{ 0<\gamma_1, \ldots , \gamma_s \leq 1} \; | 
  \Delta((\bx_{n})_{n=1}^{N} , J_{\bgamma})|/N .
\end{equation}

Let $((\bx_{n})_{n \geq 1})$ be an arbitrary sequence in  $[0,1)^s$. According to the well-known conjecture (see e.g. [1, p.67], [3, p.32])
\begin{equation} \label{2}
  \overline{\lim}_{N \to \infty}
N (\ln N)^{-s} \emph{D}^{*}((\bx_n)_{n=0}^{N-1})  >0 .
\end{equation}

{\bf Definition 1}. {\it Let  $b \geq 2, s \geq 1$, and   $0 \leq u \leq m$  be  integers 
 and let $\be = (e_1, \cdots, e_s) \in \NN^s$.
A $(u,m, \be, s)$-{\sf net} in base $b$  is a point set $\cP$ of $b^m$ points in $[0, 1)^s$ 
 such that 
 every  subinterval $J \subseteq [0, 1)^s$  of volume
${\rm Vol}(J) \geq  b^{u-m}$ which has the form
$  J = \prod_{1 \leq i \leq s} [a_ib^{-d_i},(a_i+1)b^{-d_i}),   $
with integers 
$d_i \geq 0, \;  0  \leq a_i < b^{d_i}$ and $e_i |d_i$ for $1 \le i \le s$,
contains
exactly $b^m{\rm Vol}(J)$ points of $\cP$.} \\

 If $e = (e_1, ..., e_s) = (1, ... ,1)$, we obtain a classical $(u,m, s)$-net.
For $x =\sum_{j \geq 1}  x_{j}p_i^{-j}$, 
where $x_{i} \in  Z_b =\{0,...,b-1\}$ and $m \in \NN$,   we define the truncation $[x]_m =\sum_{1 \leq j \leq m}  x_{j}b^{-j}$. If $\bx = (x^{(1)}, . . . , x^{(s)})  \in [0, 1)^s$, then the truncation $[\bx]_m$ is defined coordinatewise, that is, $[\bx]_m = 
( [x^{(1)}]_m, . . . , [x^{(s)}]_m)$.\\

{\bf Definition 2}.   {\it Let  $b \geq 2, s \geq 1$, and   $0 \leq u \leq m$  be  integers 
 and let $\be = (e_1, \cdots, e_s) \in \NN^s$.
A sequence $\bx_0, \bx_1,...$ of points in $[0, 1)^s$  is a $(u, \be, s)$-{\sf sequence} in
base $b$ if for all integers $k \geq0$ and $m > u$ the points $[\bx_n]_{m}$ with
$kb^m  \leq n < (k + 1)b^m$ form a $(u,m, \be, s)$-net in base $b$.} \\

If $e = (e_1, ..., e_s) = (1, \cdots ,1)$, we obtain a classical $(u, s)$-sequence. 
 For $x =\sum_{j \geq 1}  x_{i}p_i^{-i}$, and $\gamma =\sum_{j \geq 1}  \gamma_{i}p_i^{-i}$ 
where $x_{i},\gamma_i \in  Z_b$, we define the ($b$-adic) digitally shifted point $v$ by
$v = x \oplus \gamma := \sum_{j \geq 1}  v_{i}p_i^{-i}$,
 where $v_i \equiv x_i + \gamma_i (\mod b)$ and $v_i \in Z_b$.
For higher dimensions $s > 1$ let $\bgamma = (\gamma^{(1)}, . . . , \gamma^{(s)}) \in [0, 1)^s$. For $\bx = (x^{(1)}, . . . , x^{(s)}) \in [0, 1)^s$ we define the ($b$-adic) digital shifted point $\bv$ by 
$ \bv =\bx \oplus \bgamma =(x^{(1)} \oplus \gamma^{(1)}, . . . ,x^{(s)} \oplus \gamma^{(s)})$.
 For $n_1,n_2 \in [0,b^m)$, we define 
$n_1 \oplus n_2 := (n_1 /b^m\oplus n_2/b^m)b^m$. \\

For $x =\sum_{j \geq 1}  x_{i}p_i^{-i}$, 
where $x_{i} \in  Z_b$,  $x_i=0 $ $(i=1,...,k)$ and $x_{k+1} \neq 0$ we define the
absolute valuation  $\left\|.  \right\|_b $ of $x$ by  $\left\|x  \right\|_b =b^{-k-1}$.
Let $\left\| n  \right\|_b =b^k$ for $n \in [b^k,b^{k+1})$.\\

	 {\bf Definition 3.} {\it A digital point set   $ (\bx_{n})_{0 \leq n <b^m} $ 
	in $[0,1)^s$ is  $d-${\sf admissible} in
base $b$  if}
\begin{equation}  \label{3}
  \min_{0 \leq k <n < b^m} \left\| \bx_n \ominus \bx_k  \right\|_b
  > b^{-m-d}  \quad {\rm where} \quad \left\| \bx  \right\|_b := \prod_{i=1}^s 
	\left\|x^{(i)}  \right\|_b . 
\end{equation}
{\it A sequence  $ (\bx_{n})_{n \geq 0} $
	in $[0,1)^s$ is  $d-${\sf admissible} in
base $b$ if} $ \inf_{n >k \geq 0}
	\left\| n \ominus k
	\right\|_b \left\| \bx_n \ominus \bx_k  \right\|_b  \geq b^{-d}$. \\

By [3, p. 60] $ N  \emph{D}^{*}((\beta_n)_{n=0}^{N-1}) =O((\ln N)^{s})$ for every    $(t,s)$-sequence $(\beta_n)_{n \geq 0}$.
  In this paper we prove that   this estimate is exact for digitally shifted
	$d-$admissible $(t,s)$ sequences and in particulary for digitally shifted Niederreiter's sequence (see e.g. [1]-[4]). This result supports the conjecture
	(\ref{2}). In [2] we prove that $(t,s)$ sequences from  [1, Section 8]
are $d-$admissibles. \\

{\bf Theorem 1.} {\it  Let $s \geq 2$,  $d \geq 1$, $E_m =\{ [y]_m \;| \;y \in [0,1)\}$, 
$ (\bx_n)_{0 \leq n <b^m} $  be a $d-$admissible   $(t,m,s)$ net in base $b$, 
 $m \geq 9(d+t)(s-1)^2$. Then}
\begin{equation}\nonumber
      \max_{ \bw \in E_m^s}   b^m \emph{D}^{*}((\bx_n \oplus \bw)_{0 \leq  n < b^m}) \geq   
		 b^{-d}K_{d,t,s}^{-s+1} m^{s-1} \quad 
		\with \quad K_{d,t,s} =  4 (d+t)(s-1)^2.
\end{equation} \\

{\bf Theorem 2.} {\it Let $s \geq 1$, $d \geq 1$, $ (\bx_n)_{ n \geq 0} $  be an $d-$admissible   $(t,s)$ sequence in base $b$. Then }
\begin{equation}  \label{10}
    1+  \min_{0 \leq Q <b^m} \max_{1 \leq N \leq b^m, \bw \in E_m^{s}} 
		N \emph{D}^{*}((\bx_{n\oplus Q} \oplus \bw)_{0 \leq  n < N}) \geq   
					 b^{-d}K_{d,t,s+1}^{-s} m^{s} \quad \for \quad  m \geq 9(d+t)s^2.
\end{equation} \\

{\bf Theorem 3.} {\it Let $s \geq 1$, $ (\bx_n)_{ n \geq 0} $ be a generalized Niederreiter sequence with generating polynomials $p_1,...,p_s$, (see [1, p.266], [4, p. 242]),
$\e_i=\deg(p_i)$ $1 \leq i \leq s$,
  $e_0=e_1+ \cdots +e_s$, $d=e_0$, $t=e_0 -s$.  Then (\ref{10}) hold.} 

%
%
\newpage
\section{Proof } 
{\bf Lemma 1.} {\it Let ${\dot{s}} \geq 2$, $d \geq 1$, $ (\bx_n)_{ 0 \leq n < b^m} $  be an $d-$admissible  $(t,m,{\dot{s}})$ net in base $b$, 
  $d_0 =d+t$, $\hat{e} \in\NN$,  $0<\epsilon \leq (2d_0 \hat{e} (\dot{s}-1))^{-1}$, $\dot{m} = [m \epsilon]$,
 $\ddot{m}_i=0$,
 $\dot{m}_{i} = d_0\hat{e}\dot{m}$ $(1 \leq i \leq {\dot{s}}-1)$, 
$\ddot{m}_{\dot{s}} =m- ({\dot{s}}-1)\dot{m}_1 -t\geq 1 $,
 $\dot{m}_{\dot{s}} = \ddot{m}_{\dot{s}} +\dot{m}_1$, 
$B_{i} \subset \{0,...,\dot{m}-1\} $
$(1 \leq i \leq {\dot{s}}) $,  $\bw \in E_m^{{\dot{s}}}$
   and let
$\gamma^{(i)}=\gamma^{(i)}_1/b+...+\gamma^{(i)}_{\dot{m}_i}/b^{\dot{m}_i}$,
\begin{equation}\label{In10}
\gamma^{(i)}_{\ddot{m}_i +d_0(\hat{j}_i\hat{e} +\breve{j}_i) +\check{j}_i} =0 \quad \for\quad  1 \leq \check{j}_i < d_0, 
\qquad \qquad
 \; \gamma^{(i)}_{\ddot{m}_i +d_0(\hat{j}_i\hat{e} +\breve{j}_i) +\check{j}_i}=1 \quad \for\quad \check{j}_i =d_0 , 
\end{equation}  
and $ \hat{j}_i  \in \{0,...,\dot{m}-1\} \setminus B_{i}, \; 
  0 \leq \breve{j}_i < \hat{e}, \; 1 \leq i \leq \dot{s} $, 
 $\bgamma =(\gamma^{(1)},...,\gamma^{(\dot{s})})$, 
$B =\#B_1+...+\#B_{\dot{s}}$.
 Let there exists $n_0 \in [0,b^m)$ such that 
$ [(\bx_{n_0} \oplus \bw)^{(i)}]_{\dot{m}_i} =\gamma^{(i)}$, $1 \leq i \leq \dot{s}$, and
$m \geq 4 \epsilon^{-1}(\dot{s} -1)(1+\dot{s} B) +2t$. Then }\\
\begin{equation}  \nonumber
\tilde{\Delta}:= \Delta((\bx_n \oplus \bw)_{0 \leq n < b^m}, J_{\bgamma})  \leq
			 		-	b^{-d} \big( \hat{e}\epsilon(2({\dot{s}}-1))^{-1} \big)^{{\dot{s}}-1} m^{{\dot{s}}-1} 
			+b^{t+s}d_0\hat{e}  Bm^{\dot{s}-2}.
\end{equation} \\

{\bf Proof.} 
 Let $\br=(r_1,...,r_{\dot{s}}) \in \NN^{\dot{s}}$, $r_0=r_1+ \cdots +r_{\dot{s}}$,
$ A = \{ \br \; | \; 1 \leq r_i \leq \dot{m}_i, \; i=1,...,\dot{s} \; {\rm and}\; 
 \gamma^{(1)}_{r_1}...\gamma^{(\dot{s})}_{r_{\dot{s}}} \neq 0  \} $,  
$\dot{A}=\{ \br \in A  | \; \exists i \in [1,\dot{s}]:\; [(r_i-\ddot{m}_i-1)/(d_0\hat{e})] \in B_i  \}$, 
$A_1 =\{ \br \in A \; | \; r_0 \leq m-t\}$,
 $A_2 =\{ \br \in A \cap \dot{A} \; | \; r_0 > m-t\}$, 
 $A_3 =\{ \br \in A \setminus \dot{A} \; | \;m-t < r_0 < m+d\}$ and
$A_4 =\{ \br \in A \setminus \dot{A}\; | \; r_0 \geq m+d\}$.
We have $A =A_1 \cup A_2 \cup A_3 \cup A_4 $.
Let  
\begin{equation}  \nonumber
J_{\bgamma} =\prod_{1\leq j \leq \dot{s}} [0, \gamma^{(i)})  \qquad  \and  \qquad J_{\br,\bgamma,\bg} =\prod_{1\leq j \leq \dot{s}}
\big[ [\gamma^{(i)}]_{r_i-1}  +g_ib^{-r_i}, [\gamma^{(i)}]_{r_i-1}  +(g_i+1)b^{-r_i}   \big) .
\end{equation}
  Similarly to [3, p. 37,38], from (\ref{1}) we have that 
	$\tilde{\Delta}  =\Delta_1+\Delta_2 + \Delta_3+\Delta_4$, where \\
\begin{equation}\nonumber
 \Delta_i = \sum_{\br \in A_i} \Psi_{\bgamma}^{(\br)} ,  \qquad
    \Psi_{\br,\bgamma}  =\sum_{0 \leq g_i < \gamma^{(i)}_{r_i},\;1 \leq i \leq \dot{s}} \Psi_{\br,\bgamma,\bg}  \quad {\rm and}  \quad \Psi_{\br,\bgamma,\bg}
	= \sum_{0 \leq n < b^m} \big({\bf 1} (\bx_n  \oplus \bw , J_{\br,\bgamma,\bg} )  -b^{-r_0}\big). 
\end{equation} 
 Consider $\Delta_1$. Bearing in mind that $(\bx_n \oplus \bw)_{0 \leq n < b^m}$ is a $(t,m,{\dot{s}})$ net, we obtain $\Psi_{\br,\bgamma,\bg}  =0$. Hence  $\Delta_1=0$. \\
Consider $\Delta_2$. It is easy to verify that $\Delta_2 \leq b^{t+\dot{s}-1} d_0\hat{e}Bm^{{\dot{s}}-2}$.\\
Consider $\Delta_3$.
 We see that $r_0  \in (m-t, m+d)$. Hence $r_{\dot{s}}=r_0-r_1- ...-r_{{\dot{s}}-1} 
>m-t-(\dot{s}-1)\dot{m}_1 =\ddot{m}_{\dot{s}}$.  
 Taking into account that  $\gamma^{(i)}_{r_i}\neq 0$ and $ [(r_i -\ddot{m}_{i}-1)/(d_0\hat{e})] \notin B_i $, 
we get   $r_i=\ddot{m}_{i}+d_0j_i$ 
with some $ j_i \geq 1 $,  $1 \leq i \leq {\dot{s}}$.  Hence 
\begin{equation}\nonumber
r_0 = r_1+ \cdots +r_{\dot{s}}= \ddot{m}_{\dot{s}} +d_0(j_1 + \cdots +j_{\dot{s}}) =
 m-t + d_0(j_1 + \cdots +j_{\dot{s}} -({\dot{s}}-1) \hat{e}\dot{m})    >m-t. 
\end{equation}
Thus $r_0 \geq   m-t+d_0 =m+d$. We have a contradiction. 
Hence $ A_3 =\emptyset$ and $\Delta_3=0$.\\
Consider $\Delta_4$. Suppose that
 ${\bf 1} ( \bx_k  \oplus \bw, J_{\br,\bgamma,\bs})  =1$ for some $k \in [0, b^m)$ and  $r_0  \geq m+d$.  Then $ [ (\bx_k \oplus \bw)^{(i)}]_{r_i}    =  [\gamma^{(i)}]_{r_i}   -b^{-r_i} $,  $i=1,...,{\dot{s}}$.
Hence  $  x^{(i)}_{k,j} \ominus x^{(i)}_{n_0,j} =0$ for $j \in [1,r_i)$, $i=1,...,{\dot{s}}$. 
Therefore
\begin{equation}\nonumber
\left\| x^{(i)}_k  \ominus x^{(i)}_{n_0}  \right\|_b  \leq b^{-r_i}\quad \for \quad  i=1,...,{\dot{s}} \quad \ad 
\quad  \left\| \bx_k \ominus \bx_{n_0}  \right\|_b \leq b^{-r_0} \leq b^{-m-d}.
\end{equation}
By  (\ref{3}) and conditions of the lemma, we have a contradiction.  
Thus ${\bf 1} (\bx_k  \oplus \bw, J_{\br,\bgamma,\bs}) =0$.

We have 
$\Delta_4 \leq  -\sum_{\br \in A_4 } b^{m-r_0}$. 
 We derive $\Delta_4 \leq  -b^{-d} \#A_5$ with
$A_5 = \{ \br \in A_4\;|\; r_0=m+d\}.$  \\
Let $\hat{j}_{i} \in \{0,..., \dot{m}-1 \} \setminus B_i$, $\breve{j}_i \in [0, \hat{e}-1]$ and $r_i= \ddot{m}_{i}+d_0(\hat{e} \hat{j}_{i} +\breve{j}_{i}+1)$ for $i \in [1,\dot{s}] $. By (\ref{In10}), we get $\gamma^{(i)}_{r_i} =1$  for $i \in [1,\dot{s}] $.
Hence $A_5 \supseteq   A_6$, where
\begin{equation}\nonumber
   A_6 =\{ \br  \; |  \; r_0=m+d, \; r_i= \ddot{m}_{i}+d_0(\hat{e} \hat{j}_{i} +\breve{j}_{i}+1), \;
		\;  \hat{j}_{i} \in \{0,..., \dot{m}-1 \} \setminus B_i, \; 
		 \breve{j}_i \in [0, \hat{e}-1], \;i \in [1,\dot{s}]  \}. 
\end{equation}
Let $j_i = \hat{e} \hat{j}_i +\breve{j}_{i}+1$, for $i \in [1,\dot{s}]$.
We have  
\begin{gather}\nonumber
   r_0= \ddot{m}_{\dot{s}} +   d_0 (j_1+...+j_{{\dot{s}}}) = 
	  m-t +d_0(j_1+...+j_{\dot{s}} -(\dot{s} -1)\hat{e}\dot{m}) =m+d
	\;   \with \; d_0 =d+t.\\
{\rm Hence} \; 
j_{\dot{s}}=(\dot{s} -1)\hat{e}\dot{m} +1- j_1-...-j_{{\dot{s}}-1}.\; {\rm It\; is\; easy \;to \; verify \;that} \;j_{\dot{s}}\in [1,\hat{e}\dot{m}]\;\for\;  \hat{j}_{i} \in [\dot{m}- \dddot{m}, \dot{m}-1],
	 \qquad \qquad\qquad \qquad\qquad \qquad\qquad \qquad\qquad \qquad\qquad \qquad \nonumber \\
\for\; i \in[1,{\dot{s}}-1],\; \with  \;\dddot{m}=[\dot{m}/({\dot{s}}-1)].\; {\rm Thus} \;\#A_6 \geq \#A_7,\; \where
\qquad \qquad\qquad \qquad\qquad \qquad\qquad \qquad \qquad\qquad \qquad
 \nonumber
\end{gather}
\begin{equation}\nonumber
   A_7 =\{ (j_1,..., j_{{\dot{s}}-1} )  \; | \;  j_i = \hat{e} \hat{j}_i +\breve{j}_{i}+1, \quad
			 \hat{j}_{i} \in \{0,..., \dot{m}-1 \} \setminus B_i,\quad \breve{j}_i \in [0, \hat{e}-1],
			\quad	i \in [1,\dot{s}], 
\end{equation}
\begin{equation}\nonumber
  \hat{j}_{i} \in [\dot{m}- \dddot{m}, \dot{m}-1],\quad i \in [1,\dot{s}-1] 
	\qquad \ad \qquad
	j_{\dot{s}}=(\dot{s} -1)\hat{e}\dot{m} +1- j_1-...-j_{{\dot{s}}-1}
		  \}.  
\end{equation}
We obtain $\#A_7 \geq \#A_8 - \hat{e}\#B_{\dot{s}}m^{{\dot{s}}-2}$, where
\begin{equation}\nonumber
   A_8 =\{ (j_1,..., j_{{\dot{s}}-1} )   |   j_i = \hat{e} \hat{j}_i +\breve{j}_{i}+1, \;
			 \hat{j}_{i} \in \{\dot{m}- \dddot{m},..., \dot{m}-1 \} \setminus B_i, \breve{j}_i \in [0, \hat{e}-1] ,\; i \in [1,\dot{s}-1] \}.
\end{equation}
Therefore
\begin{equation}\nonumber
    \# A_8 \hat{e}^{{-\dot{s}}+1} \geq \# \{ (\hat{j}_1,..., \hat{j}_{{\dot{s}}-1} )  \; | \;  
		  1 \leq \hat{j}_i  \leq  \dddot{m} - \#B_i,\; 1 \leq i \leq{\dot{s}}-1 \}  
			\geq  (\dddot{m}-B)^{{\dot{s}}-1} 
\end{equation}
\begin{equation}\nonumber
    = \dddot{m}^{{\dot{s}}-1} (1-B/\dddot{m})^{{\dot{s}}-1} \geq 
			\dddot{m}^{{\dot{s}}-1} \big(1-({\dot{s}}-1)B/\dddot{m}\big) \geq   
			 \big( m \epsilon(2({\dot{s}}-1))^{-1} \big)^{{\dot{s}}-1}
			 			- (\dot{s} -1) B \dddot{m}^{{\dot{s}}-2} 	
\end{equation} 
for  $m \geq 4 \epsilon^{-1}(\dot{s} -1)(1+\dot{s} B) +2t$.
Therefore  $\tilde{\Delta} \leq -b^{-d} \big( \hat{e}\epsilon(2({\dot{s}}-1))^{-1} \big)^{{\dot{s}}-1} m^{{\dot{s}}-1} 
			+b^{t+s} d_0
			\hat{e} B m^{\dot{s}-2} $.\\
Thus Lemma 1 is proved. \qed  \\

{\bf Proof of Theorem 1.}  Using Lemma 1 with ${\dot{s}}=s $, $B_{i} =\emptyset$ $(1 \leq i \leq s)$, $B=0$,
  $\hat{e} =1$, $\epsilon= (2(s-1)d_0)^{-1} $, 
$n_0 = 0$, and $\bw =[\bgamma \ominus \bx_0]_m$, 
    we obtain the assertion of Theorem 1.           \qed 

{\bf Proof of Theorem 2.} 
  According to  [3, Lemma 3.7], we have
\begin{equation}  \nonumber
 1+  \sup_{1 \leq N \leq b^m}  N \emph{D}^{*}((\bx_{n\oplus Q} \oplus \bw )_{n=0}^{N -1})
 \geq   b^m \emph{D}^{*}((  \bx_{n\oplus Q} \oplus \bw, n/b^m )_{n=0}^{b^m -1}) 
		=   
		b^m \emph{D}^{*}((  \bx_n\oplus \bw, (n \ominus Q)/b^m )_{n=0}^{b^m -1}).
\end{equation} 
By  (\ref{3}) and [1, Lemma 4.38], we have that $((  \bx_n, n/b^m )_{0 \leq  n < b^m})$   is a $d-$admissible $(t,m,s+1)-$net in base~$b$.
Using Lemma 1
 with ${\dot{s}}=s+1 $,   $x_n^{(s+1)} =n/b^{m}$, $B_{i} =\emptyset$ $(1 \leq i \leq s+1)$, $B=0$,
  $\hat{e} =1$, $\epsilon= (2sd_0)^{-1} $, 
$n_0 = Q \oplus \gamma^{(s+1)}b^m$,  and $\bw =([(\gamma^{(1)},...,\gamma^{(s)}) \ominus \bx_{n_0}]_m,{-Q/b^m})$,
 we get the assertion of Theorem~2.           \qed

{\bf Lemma 2.} {\it Let $(\bx_n)_{n \geq 0}$ be a $(0, \be,s)$ sequence in base $b$. 
Then $(\bx_n)_{n \geq 0}$ is $e_0-$admissible.} \\

{\bf Proof.} Suppose that $(\bx_n)_{n \geq 0}$  is not a $e_0$ admissible. Then there exists
 $ n_0 > k_0 \geq 0$ with 
$\left\| n_0 \ominus k_0 \right\|_b $  $\times  \left\| \bx_{n_0} \ominus \bx_{k_0} \right\|_b  \leq 
b^{-e_0-1}$. Let   $\left\| n_0 \ominus k_0 \right\|_b= b^{\tilde{d}}$, and let
 $\left\| x_{n_0}^{(i)} \ominus x_{k_0}^{(i)} \right\|_b =b^{-d_i-1}$ $ (i=1,...,s)$.
Hence $\varkappa := \tilde{d}  - \sum_{1 \leq i \leq s} 
			( d_i +1) +e_0 +1 \leq 0$. 
Let $\dot{d}_i =[d_i/e_i]e_i \geq d_i -e_i +1 $,
  $a_i =[x_{n_0}^{(i)}]_{\dot{d}_i} b^{\dot{d}_i}$ $(i=1,...,s)$  and let 
$J =\prod_{1 \leq i \leq s} [a_i b^{-\dot{d}_i}, (a_i +1)b^{-\dot{d}_i}) $.
We have $x_{n_0,j}^{(i)} =x_{k_0,j}^{(i)}$ for all $j \in [1,d_i]$, $i \in [1,s]$.
Hence $\bx_{n_0} , \bx_{k_0} \in J$. 
 We derive
\begin{equation}  \label{70}
  0 \geq  \varkappa  =  \tilde{d} +1 - \sum_{1 \leq i \leq s} 
			(d_i -e_i +1) \geq \tilde{d} +1 - \sum_{1 \leq i \leq s} \dot{d}_i,
			 \quad \ad  \quad 1 \geq b^{\tilde{d} +1}{\rm Vol}(J).
\end{equation}
 Let $n_0 = \dot{n}_0 b^{\tilde{d} +1} + \ddot{n_0} $ where 
$\ddot{n_0} \in [0, b^{\tilde{d} +1})$.  It is easy to see that 
$k_0 = \dot{n}_0 b^{\tilde{d} +1} + \ddot{k_0} $  with some 
$\ddot{k_0} \in [0, b^{\tilde{d} +1})$. Hence 
 $n_0, k_0 \in [ \dot{n}_0 b^{\tilde{d} +1}, (\dot{n}_0 +1)b^{\tilde{d} +1} ) =:W$.
Thus $ \sum_{n \in W}  {\bf 1} ( \bx_n,J)
  \geq 2$. Bearing in mind (\ref{70}), we obtain that
 $(\bx_n)_{n \geq 0}$  is not  $(0,\be,s)$ sequence. We have a contradiction.
Hence Lemma 2 is proved. \qed

{\bf Proof of Theorem 3.} Let $\be =(e_1,...,e_s)$. By [4] and [1, p.266], we have that   $ (\bx_n)_{ n \geq 0} $ is a $(0,\be, s) $ and $(e_0-s, s) $ 
 sequence. Applying Lemma 2 and Theorem 2,  we obtain the assertion of Theorem 3.           \qed 

\end{document}